\documentclass[12pt,a4]{amsart}
\usepackage{amssymb,mathrsfs,amsmath,dsfont,array}

\usepackage{multirow,arydshln}
\usepackage{graphicx}
\usepackage{bm}
\usepackage{tikz}
\usepackage{capt-of}
\usepackage{mathbbol}     
\usepackage{color}
\usepackage{verbatim} 
\usepackage[all,cmtip]{xy}
\usepackage{amsbsy}
\usepackage{amsmath}
\usepackage{nccmath}

\usepackage{hyperref}

\newtheorem{Thm}[equation]{Theorem}
\newtheorem{Pro}[equation]{Proposition}

\newtheorem{lem}[equation]{Lemma}
\newtheorem{exa}[equation]{Example}

\newtheorem{deft}[equation]{Definition}
\newtheorem{rem}[equation]{Remark}

\numberwithin{equation}{section}\makeatletter
\newcommand{\AlignFootnote}[1]{%
    \ifmeasuring@
    \else
        \footnote{#1}%
    \fi
}
\makeatother
\numberwithin{equation}{section}

\newcommand{\End}{\mathop{{\rm End}}\nolimits}

\renewcommand{\hat}{\widehat}

\def\top{\hbox{\tiny$\op$ }}
\def\bs{\boldsymbol}

\def\a1s{{\dot A},\cdots, a_s}
\def\a{\alpha}

\def\andd{\quad\hbox{and}\quad}

\def\b{\beta}

\def\bl4{B_{\ell\geq4}}

\def\bbbc{{\mathbb C}}

\def\d{\delta}
\def\D{\Delta}

\newcommand{\rmd}{{\tt  d}}

\def\fg{\mathfrak{g}}

\def\scg{\mathscr{G}}

\def\fh{\mathfrak{h}}
\def\sch{\mathscr{H}}

\def\ii{\mathcal{I}}

\def\fk{\mathfrak{k}}

\def\lam{\lambda}
\def\Lam{\Lambda}

\def\ep{\epsilon}

\def\1k{\frac{1}{k}}
\def\op{\oplus}
\def\ot{\otimes}

\def\sub{\subseteq}
\def\sg{\sigma}

\def\pf{\noindent{\bf Proof. }}

\def\sspan{\hbox{\rm span}}

\def\ft{\mathfrak{t}}

\def\T{{\mathcal T}}

\def\bbbz{{\mathbb Z}}

\def\1il{1\leq i\leq\ell}

\renewcommand{\hat}{\widehat}
\def\bs{\boldsymbol}

\newcommand{\Bigop}[2]{\raisebox{0.2ex}{\scalebox{.7}{$\displaystyle \bigoplus_{#1}^{#2}\;$}}}

\newcommand{\tinyop}{\raisebox{0.2ex}{\scalebox{.7}{$\displaystyle ~\oplus\;$}}}
\newcommand{\smpro}{\raisebox{0.2ex}{\scalebox{.6}{$\displaystyle ~\#\;$}}}

\def\ddk{\frac{\partial}{\partial\zeta_k}}

\def\bddl{{\bs{ \partial}}_l}
\def\bddp{{\bs{ \partial}}_p}
\def\bddq{{\bs{ \partial}}_q}
\def\bddti{\bs{\rmd}_i}
\def\bddtj{\bs{\rmd}_j}
\def\bddtr{\bs{\rmd}_r}
\def\bddts{\bs{\rmd}_s}

 \title{Quasi-Poisson Modules and  Harish-Chandra  $\bs{AD}$-Modules}

\thanks{The author would like to thank Professor Yuly Billig for helpful discussion during her visit to Carleton University. Also her  gratitude to Professor Eswara Rao Senapathi for his helpful comments on the manuscript.}
\thanks{This  research was in part supported by a grant from IPM (No. 1405170411)}
\begin{document}

\maketitle

\pagestyle{myheadings}
\markboth{}{}

\centerline{ Malihe Yousofzadeh}

\centerline{{\scalebox{0.65} {Department of Pure Mathematics, Faculty of Mathematics and Statistics, University of Isfahan,}}}\centerline{{\scalebox{0.65} { P.O.Box 81746-73441, Isfahan, Iran, and }}}

\centerline{{\scalebox{0.65} { School of Mathematics,
 Institute for Research in
	Fundamental Sciences (IPM), }}}
\centerline{{\scalebox{0.65} {P.O. Box: 19395-5746, Tehran, Iran.}
}}

\centerline{{\scalebox{0.65} {		 ma.yousofzadeh@sci.ui.ac.ir \& ma.yousofzadeh@ipm.ir.	}}}

\curraddr{}


\begin{abstract}
We introduce the notion of quasi-Poisson modules over Lie-Rinehart pairs and prove that for the Lie-Rinehart pair  $(\dot A,\dot\fk)$  in which $\dot A:=\bbbc[t_1^{\pm1},\ldots,t_m^{\pm1}]\ot\Lam_n$ and $\dot\fk:={\rm Der}(\dot A)$, there is a one-to-one correspondence between simple cuspidal quasi-Poisson modules over $(\dot A,\dot\fk)$   and simple  cuspidal Harish-Chandra $A\fk$-modules  for $A:=\bbbc[t_0^{\pm1}]\ot \dot A$ and $\fk:={\rm Der}(A).$ We also classify simple cuspidal quasi-Poisson modules over the Lie-Rinehart pair  $(\dot A,\dot\fk)$ and show that each such  module is a tensor module $\dot A\ot \Omega$ for  an admissible $\frak{gl}(m+1,n)$-module $\Omega$ via a prescribed  action. In particular, we get that every nontrivial simple cuspidal  Harish-Chandra
$\fk$-module  is a $\fk$-quotient of a tensor $A\fk$-module. It generalizes the known result that  every nontrivial simple cuspidal strong  Harish-Chandra
$\fk$-module  is a $\fk$-quotient of a tensor $A\fk$-module.
\end{abstract}
\maketitle

\section{Introduction}
The kernel of the universal central extension of $\scg\ot \mathscr{A}$  where $\scg$ is  a finite dimensional simple Lie algebra  and $\mathscr{A}$ is  a unital commutative associative algebra, is the K\"ahler differentials  $\Omega(\mathscr{A})$ of $\mathscr{A}$ modulo the exact K\"ahler differentials $d(\mathscr{A})$; see  \cite{KL, NY}. A toroidal Lie algebra is  the universal central extension $(\scg\ot\dot{\mathscr{A}})\op\frac{\Omega(\dot{\mathscr{A}})}{d({\dot{\mathscr{A}}})}$ of $\scg\ot \dot{\mathscr{A}}$ where  $\dot{\mathscr{A}}$ is the set of  Laurent polynomials in $\ell$  variables $t_1,\ldots, t_\ell$.
The split  central extension ${\rm Der}(\dot{\mathscr{A}})\op\frac{\Omega(\dot{\mathscr{A}})}{d({\dot{\mathscr{A}}})}$ of $\dot{\mathscr{D}}:={\rm Der}(\dot{\mathscr{A}})$ appears in generalizing the vertex operator construction on Fock spaces \cite{EM}. 
 Then, one can construct 
 \[\tau:=(\fg\ot\dot{\mathscr{A}})\op \dot{\mathscr{D}}\op\frac{\Omega(\dot{\mathscr{A}})}{d({\dot{\mathscr{A}}})}\]
 which is known in the literature  as a full toroidal Lie algebra. In \cite{JM}, the authors prove that the classification of simple modules over $\tau$ is reduced to the classification of   simple so-called $\dot{\mathscr{A}}\dot{\mathscr{D}}$-modules. In \cite{Lar}, the author constructs a class of simple $\dot{\mathscr{A}}\dot{\mathscr{D}}$-modules referred to as Shen-Larsson modules and in \cite{E1}, the author shows that Shen-Larsson modules exhaust all   simple $\dot{\mathscr{A}}\dot{\mathscr{D}}$-modules having finite weight space decomposition with respect to $\{\rmd_i:=t_i\frac{d}{dt_i}\mid 1\leq i\leq \ell\}$ with bounded weight multiplicities. 
The pair $(\dot{\mathscr{A}},\dot{\mathscr{D}})$ is an example of a  Lie-Rinehart pair $(A,\fk)$
 which contains  a unital supercommutative associative superalgebra $(A,\cdot)$ as well as  a Lie superalgebra $(\fk,[\cdot,\cdot])$  satisfying the following:
\begin{itemize}
\item $\fk$ is an $A$-module; we denote the action of $a\in A$ on  $\d\in \fk$   by $a\d,$
\item $\fk$ acts on $A$ as derivations;  we denote the action of $\d\in \fk$ on $a\in A$ by $\d(a),$
\item for $a\in A$ and $\d,\gamma\in \fk,$ we have $[a\d,b\gamma]=a\d(b)\gamma-(-1)^{|a\d||b\gamma|}b\gamma(a)\d+(-1)^{|b||\d|}ab[\d,\gamma].$
\end{itemize}

For a   Lie-Rienhart pair $(A,\fk)$, we say that a superspace $\Omega$ is an $A\fk$-module if $\Omega$ is both an $A$-module and  a $\fk$-module satisfying the Leibniz rule
\[\d(a\omega)=\d(a)\omega+(-1)^{|\d||a|}a(\d \omega)\quad(a\in A,~\d\in \fk,~\omega\in\Omega).\] 
If $\Phi:U(\fk)\longrightarrow U(\fk)\ot U(\fk)$ is the natural co-product of the Hopf algebra $U(\fk)$ with $\Phi(u)=\displaystyle{\sum_i u_i^{1}\ot u_i^2},$ one can construct  the smash product $A\smpro U(\fk)$ whose underlying vector space is $A\ot U(\fk)$ with the product
\[(a\ot u)\smpro (b\ot v)=\sum_i (-1)^{|b||u_i^2|}au_i^1(b)\ot u_i^2v;\] here $u_i^1(b)$ is the induced action of $U(\fk)$ on $A$. $A\fk$-modules are in fact modules over the associative superalgebra $A\smpro U(\fk)$ and they are recovered by the action of $A\smpro\fk.$

Assume $n$ is a positive integer and $\Lam_n$ is the Gra\ss{}mann algebra in $n$-variables  $\zeta_1,\ldots,\zeta_n$. Also, assume that $\mathfrak{G}$ is one of the unital associative supercommutative  superalgebras \[A:=\bbbc[t_0^{\pm1},\ldots,t_m^{\pm1}]\ot \Lam_{n}\hbox{ or }A^+:=\bbbc[t_0,\ldots,t_m]\ot \Lam_{n}\] and $$\fk_\mathfrak{G}:={\rm Der}(\mathfrak{G}).$$  
A $\fk_\mathfrak{G}$-module $V$ is called a   Harish-Chandra module (resp. a strong Harish-Chandra module) if $V$ has a weight space decomposition with respect to 
$\sch:=\{t_i\frac{d}{dt_i},\zeta_k\frac{\partial}{\partial_k}\mid 0\leq i\leq m,1\leq k\leq n\}$ (resp. $\sch':=\{t_i\frac{d}{dt_i}\mid 0\leq i\leq m\}$) with finite dimensional weight spaces.
A (strong) Harish-Chandra module is called  cuspidal   if its weight spaces have bounded weight multiplicities.

A $\mathfrak{G}\fk_\mathfrak{G}$-module $V$ is called a    (strong) Harish-Chandra module (resp. cuspidal (strong) Harish-Chandra module) if $V$ is a (strong) Harish-Chandra module (resp. cuspidal (strong) Harish-Chandra module) as a $\fk_\mathfrak{G}$-module.

Simple cuspidal strong Harish-Chandra $A\fk_A$-modules  and simple cuspidal strong Harish-Chandra $\fk_A$-modules have been studied and classified in \cite{XL1,BFIK}. Also simple cuspidal Harish-Chandra $\fk_{A^+}$-modules  and simple cuspidal strong Harish-Chandra $\fk_{A^+}$-modules  have been classified  in \cite{LX,CLX}.
The only remaining case is the classification of  simple cuspidal  Harish-Chandra $A\fk_A$-modules and $\fk_A$-modules.
In  \cite{XW}, the authors  deal with  $\fk_A$-modules; they introduce  tensor modules and investigate their simple quotients.  In this work,   we complete the classification of simple cuspidal Harish-Chandra $A\fk_A$-modules and show that every nontrivial simple cuspidal  Harish-Chandra $\fk_A$-module  is a $\fk_A$-quotient of a tensor $A\fk_A$-module.
\smallskip

We introduce the concept of a (cuspidal) QP-module over a Lie-Rinehart pair and  for
\[A=\bbbc[t_0^{\pm1},\ldots,t_m^{\pm1}]\ot \Lam_n\quad\quad\hbox{where $m,n$ are positive inetegers},\] and 
\[\dot A:=\bbbc[t_1^{\pm1},\ldots,t_m^{\pm1}]\ot \Lam_n\andd \dot\fk:={\rm Der}(\dot A),\]
we show that there is a one-to-one correspondence between simple   cuspidal Harish-Chandra $A\fk_A$-modules and simple  cuspidal QP-modules over $(\dot A,\dot\fk).$
We classify simple  cuspidal QP-modules  over $(\dot A,\dot\fk)$ and as a byproduct, we obtain the classification of simple  cuspidal  Harish-Chandra $A\fk_A$-modules. 

In the first place, we associate to a QP-module  over $(\dot A,\dot \fk)$, an $A\fk_A$-module. A simple cuspidal QP-module  over $(\dot A,\dot\fk)$ is associated with  a  simple cuspidal Harish-Chandra  $A\fk_A$-module $V$; in particular, $V$ admits  a weight space decomposition $V=\Bigop{\mu}{}V^\mu$ with respect to $\{\rmd_i=t_i\frac{d}{dt_i}\mid 0\leq i\leq m\}.$ We fix a weight space $V^\mu$ and  note that unlike  the case of strong Harish-Chandra modules, $V^\mu$ is not necessarily finite dimensional.  It is a module over the centralizer 
$\T$ of  \[A\cup    \{t_i\frac{d}{dt_i},\frac{\partial}{\partial\zeta_k}\mid 0\leq i\leq m,1\leq k\leq n\}\]
in the Lie subsuperalgebra $A\smpro\fk_A$ of ${\rm Lie}(A\smpro U(\fk_A))$.
The subalgebra $\T$ has a  filtration 
\(\ldots\sub\T^2\sub\T\)
and \[\T/\T^2\simeq\frak{gl}(m+1,n). \] 
We show that the $\T$-module $V^\mu$ contains a simple submodule $\Omega$ annihilated by $\T^2$ and so it is a $\frak{gl}(m+1,n)$-module. Moreover, we show that the module $V$ can be  recovered from $\frak{gl}(m+1,n)$-module $\Omega.$  In fact, the classification of simple cuspidal QP-modules of $(\dot A,\dot\fk)$ (equivalently, simple cuspidal Harish-Chandra $A\fk_A$-modules)  is reduced to  the classification of admissible  $\frak{gl}(m+1,n)$-modules, i.e.,
 simple modules admitting weight space decomposition with respect to the diagonal matrices with finite dimensional weight spaces and bounded weight multiplicities. We point out  that in \cite{CM}, the authors  use Kac functor to reduce  the classification of admissible $\frak{gl}(m+1,n)$-modules to the classification of admissible modules over the even part of $\frak{gl}(m+1,n)$ which is a  finite dimensional reductive Lie algebra. The latter modules have been extensively   studied in the literature; see \cite{F, Mat}

Throughout the paper, all vector spaces  and tensor products are defined on the field  $\bbbc$ of  complex numbers. The dual space of a vector space $V$ is denoted by $V^*$. For a superspace $V$,  we set $L(V):=\bbbc[t_0^{\pm1}]\ot V$ and  identify $V$ with a subset of $L(V)$ via $u\mapsto 1\ot u$.
We also  fix positive integers $m$ and $n$ {and  denote by $E_{i,j}$ ($0\leq i,j\leq m+n$) the elementary matrices. Moreover, by  $\bar i$ ($1\leq i\leq n$), we mean   $i+m$.} We set \[t^{\bs{\bar r}}:=t_0^{r_0}\cdots t_m^{r_m}\quad \hbox{for $\bs{\bar r}=(r_0,\ldots,r_m)\in\bbbz^{m+1},$}\] 
\[{t^{\bs{\bar r}}:=t_1^{r_1}\cdots t_m^{r_m}\quad \hbox{for $\bs{\bar r}=(r_1,\ldots,r_m)\in\bbbz^{m}$}}\] 
and for  a subset $I\sub\{1,\ldots,n\},$ by $\zeta_I\in \Lam_n$,  we mean   
\[\small{
\zeta_I:=\left\{\begin{array}{ll}
\zeta_{j_1}\cdots\zeta_{j_p}& \hbox{if } I=\{j_1<\cdots<j_p\},\\
1&  \hbox{if } I=\emptyset.
\end{array}\right.}\]
Finally, for two symbols ${\tt x}$ and ${\tt y},$ by 
$\d_{{\tt x},{\tt y}},$ we mean the Kronecker delta.  
\section{Certain Weight Module}\label{section1}
Suppose that $\bs{S}$ is the ideal of $A=\bbbc[t_{ 0}^{\pm1},\ldots,t_{m}^{\pm1}]\ot \Lam_n$ generated by $\{t_i-1,\zeta_k\mid {0}\leq i\leq m, 1\leq k\leq n\}$ and  set 
\begin{equation}\label{deltadelta'}
\begin{array}{rl}
\D:=&\{\bddti:=t_i\frac{d}{dt_i},\partial_{k}:=\frac{\partial}{\partial\zeta_k}\mid 0\leq i\leq m, 1\leq k\leq n\}\hbox{ as well as }\\
 \D':=&\{\frac{d}{dt_i},{\partial_k}\mid 0\leq i\leq m, 1\leq k\leq n\}.
 \end{array}
\end{equation}
The subset of $\bs{S}$ consisting of elements of the form  
\[ (t_{0}-1)^{r_1}\cdots(t_m-1)^{r_m}(t^{-1}_{0}-1)^{s_1}\cdots(t^{-1}_m-1)^{s_m}\zeta_I\]
for $r_0,\ldots r_m,s_1,\ldots,s_m \in\bbbz^{\geq0}$ and $I\sub\{1,\ldots,n\}$ with $|I|+r_0+\cdots +r_m+s_0+s_1+\cdots+s_m>0, $ is a basis for $\bs{S}.$
This gives  the filtration  $\cdots\sub \bs{S}^3\sub\bs{S}^2\sub\bs{S}^1=\bs{S}$ of ideals of $\bs{S}$ where for a positive integer $\ell,$
\begin{equation}\label{mpower}
\small{\bs{S}^\ell :=\sspan_\bbbc\{(t_{0}-1)^{r_0}\cdots(t_m-1)^{r_m}(t^{-1}_{0}-1)^{s_0}\cdots(t^{-1}_m-1)^{s_m}\zeta_I\mid \sum_{i=0}^m(r_i+s_i)+|I|\geq \ell, r_i,s_i\geq 0\}}.
\end{equation} 
We have 
\begin{equation}\label{mods2}
t^{\bs{\bar{r}}}-1\stackrel{({\rm mod}\bs{S}^2)}{\equiv\joinrel\equiv}r_0(t_{0}-1)+\cdots +r_m(t_{m}-1)\andd
t^{\bs\bar{r}}\zeta_k\stackrel{({\rm mod}\bs{S}^2)}{\equiv\joinrel\equiv}\zeta_k,
\end{equation}
for $\bs{\bar r}=(r_0,\ldots,r_m)\in\bbbz^{m+1}$ and $1\leq k\leq n.$
Next suppose  $\bs{S}^+$ is the linear span of the elements of the form
\[ (t_{0}-1)^{r_0}\cdots(t_{m}-1)^{r_m}\zeta_I\]
for $r_0,\ldots r_m\in\bbbz^{\geq 0}$ and $I\sub\{1,\ldots,n\}$ with $|I|+r_0+\cdots +r_m>0. $ Therefore, as above, we have a natural filtration $\cdots\sub (\bs{S}^+)^3\sub(\bs{S}^+)^2\sub(\bs{S}^+)^1=\bs{S}^+$ of $\bs{S}^+$ where 
\[(\bs{S}^+)^\ell :=\sspan_\bbbc\{(t_{0}-1)^{r_0}\cdots(t_m-1)^{r_m}\zeta_I\mid  \sum_{i=0}^mr_i+|I|\geq \ell, r_i\in\bbbz^{\geq 0}\}\quad(\ell\in\bbbz^{>0}).\] 

We note  that  $\bs{S}\D$, $\bs{S}\D',$ $\bs{S}^+\D$ and  $\bs{S}^+\D'$ are Lie superalgebras with respect to  the commutator bracket and that $\bs{S}^\ell\D\sub \bs{S}\D$, $\bs{S}^\ell\D'\sub \bs{S}\D'$, $(\bs{S}^+)^\ell\D\sub \bs{S}^+\D$ and  $(\bs{S}^+)^\ell\D'\sub \bs{S}^+\D'$ ($\ell\in\bbbz^{>0}$) are ideals. 
\begin{lem}\label{power-less}
Suppose that  $k<k'$ are positive integers, then,  we have 
\item[(i)]
$\bs{S}^k=(\bs{S}^+)^k+\bs{S}^{k'}$ and    $\frac{\bs{S}^k\D'}{\bs{S}^{k'}\D'}\simeq \frac{(\bs{S}^+)^k\D'}{(\bs{S}^+)^{k'}\D'} .$
\item[(ii)]
   $\bs{S}^k\D'=\bs{S}^k\D .$
\end{lem}
\pf  (i) The first assertion follows from the fact that for a positive integer $r$, 
\begin{align*}
(t_{j}^{-1}-1)^{r}\in&(1-t_{j})^{r}+{\bs{S}}^{r+1}.
\end{align*}
For the second  assertion,  we see that 
\begin{align*}
\frac{{\bs{S}}^k\D'}{{\bs{S}}^{k'}\D'}=\frac{({\bs{S}}^+)^k\D'+{\bs{S}}^{k'}\D'}{{\bs{S}}^{k'}\D'}\simeq \frac{({\bs{S}}^+)^k\D'}{{\bs{S}}^{k'}\D'\cap ({\bs{S}}^+)^k\D'}=\frac{({\bs{S}}^+)^k\D'}{({\bs{S}}^+)^{k'}\D'}.
\end{align*}
(ii) easily  follows as for $0\leq i\leq m,$ we have  
\[\hspace{1.7cm}t_i\frac{d}{dt_i}=(t_i-1)\frac{d}{dt_i}+\frac{d}{dt_i}\andd 
\frac{d}{dt_i}=t_i^{-1}t_i\frac{d}{dt_i}=(t_i^{-1}-1)t_i\frac{d}{dt_i}+t_i\frac{d}{dt_i}.\hspace{1.7cm}\qed\]

\begin{lem}\label{power-1} 
Let $k$ be a positive integer and $\ii$ be the ideal of $\bs{S}\D'$ generated by $(\bs{S}^+)^k\D'.$ We have $\bs{S}^{k+2}\D'\sub \ii.$
In particular, for   an  $\bs{S}\D'$-module $V$ with 
 $(\bs{S}^+)^k\D' V=\{0\},$  we have $\bs{S}^{k+2}\D' V=\{0\}.$
\end{lem}
\pf 
Set {\small \[D:=\sum_{i=0}^m(t_i-1)\frac{d}{dt_i}+\sum_{k=1}^n\zeta_k\partial_k\andd\] \[(\bs{S}^+)_r:=\{(t_{0}-1)^{s_0}\cdots(t_{m}-1)^{s_m}\zeta_I\mid s_0,\ldots,s_m\in\bbbz^{\geq 0},~I\sub\{1,\ldots,n\},~s_0+\cdots+s_m+|I|=r\} \]} for $r\in\bbbz^{>0}.$   The proof is verified through the following four steps (see \cite[Lem.~3.7]{XL1}):
\begin{itemize}
\item[(1)] $\bs{S}^{k+2}=(\bs{S}^+)_1\bs{S}^{k+1},$
\item[(2)] $\bs{S}^{k+1}D\sub\sspan_\bbbc\{t^{\bs{\bar r}} ZYD\mid \bs{\bar r} \in \bbbz^{m+1}, Z\in (\bs{S}^+)_1, Y\in (\bs{S}^+)^k\},$
\item[(3)] $(\bs{S}^+)_1\bs{S}^{k+1}\D'$ is a subset of the ideal generated by $\bs{S}^{k+1}D$,
\item[(4)] $\bs{S}^{k+2}\D'\sub \ii.$\hfill$\Box$
\end{itemize}

\medskip
In \cite[Lem.~3.7]{XL1},  it is shown  that for a finite dimensional simple $\bs{S}\D$-module $V$, one has  $\bs{S}^2\D V=\{0\}.$ The following proposition and theorem extend  this result to finite weight modules.
To state them  and in order to simplify the notations, we  first introduce some conventions:
For $\bs{\bar s}=(s_0,\ldots,s_m)\in(\bbbz^{\geq0})^{m+1}$ and $0\leq i\neq j\leq m,$ we set

\begin{align*}\widehat{\hspace{3mm}(t-1)_{j}^{\bs{\bar s}}}:=&(t_0-1)^{s_0}\cdots(t_m-1)^{s_m}/(t_j-1)^{s_j}\andd \\
\widehat{\hspace{3mm}(t-1)_{i,j}^{\bs{\bar s}}}:=&(t_0-1)^{s_0}\cdots(t_m-1)^{s_m}/(t_i-1)^{s_i}(t_j-1)^{s_j}.
\end{align*}
We next set 
\begin{equation}\label{h&h'}
\fh:=\sspan_\bbbc\{h_k:=\zeta_k{\partial_k}\mid 1\leq k\leq n\}\andd \fh':=\sspan_\bbbc\{h'_i:=(t_i-1)\frac{d}{dt_i}\mid 0\le i\leq m\}.
\end{equation}
 We define 
\begin{align*}
\ep_i:&\fh'\longrightarrow \bbbc,\quad
h'_j\mapsto \d_{i,j}\quad(0\leq i,j\leq m)&&\andd \d_k:&\fh\longrightarrow \bbbc,\quad
h_l\mapsto \d_{k,l}\quad(1\leq k,l\leq n).
\end{align*}
For  a subset $I\sub\{1,\ldots,n\},$ we set   
\[\small{
 \d_I:=\left\{\begin{array}{ll}
\d_{j_1}+\cdots+\d_{j_p}& \hbox{if } I=\{j_1<\cdots<j_p\}\\
0&  \hbox{if } I=\emptyset.
\end{array}
\right.}\]
We have 
\[\bs{S}^+\D'=\Bigop{\mu\in\fh^*}{ }(\bs{S}^+\D')^\mu\quad\hbox{with}\quad (\bs{S}^+\D')^\mu=\{x\in \bs{S}^+\D'\mid [h_k,x]=\mu(h_k)x\quad(1\leq k\leq n)\},\] and 
\[\bs{S}^+\D'=\Bigop{\mu\in(\fh')^*}{ }(\bs{S}^+\D')^{(\mu)}\quad\hbox{with}\quad (\bs{S}^+\D')^{(\mu)}=\{x\in \bs{S}^+\D'\mid [h'_i,x]=\mu(h'_i)x\quad(0\leq i\leq m)\}.\] 
In fact, for  $\bs{\bar s}=(s_0,\ldots,s_m) \in (\bbbz^{\geq0})^{m+1}$,  $I\sub\{1,\ldots,n\}$, $0\leq j\leq m$ and $1\leq k\leq n,$ we have  
\begin{align*}
&\small{(t_0-1)^{s_0}\cdots(t_m-1)^{s_m}\zeta_{I}{\partial_k}\in (\bs{S}^+\D')^{(s_0\ep_1+\cdots+s_m\ep_m)}\cap(\bs{S}^+\D')^{\d_I-\d_k}\andd}\\ 
&\small{(t_0-1)^{s_0}\cdots(t_m-1)^{s_m}\zeta_{I}\frac{d}{dt_j}\in (\bs{S}^+\D')^{(s_0\ep_1+\cdots+s_m\ep_m-\ep_j)}}\cap (\bs{S}^+\D')^{\d_I}.
\end{align*}
In particular, for $0\leq i,j\leq m,$ $p,p'\in\bbbz^{>0}, \bs{\bar  s}=(s_0,\ldots,s_m)\in(\bbbz^{\geq0})^{m+1}$, $1\leq k\leq n$, and $I\sub\{1,\ldots,n\}$, setting
{\footnotesize\begin{equation}\label{new-MM-1}\begin{array}{rlrl}
\hspace{3mm}\partial_{i,j,p}^{\bs{\bar s}}&:=(t_j-1)^p\widehat{\hspace{1mm}(t-1)_j^{\bs{\bar s}}}\frac{d}{dt_i},&\partial_{j,p,k,I}&:=(t_j-1)^p\zeta_I{\partial_k},\\\\
\hspace{3mm}\partial_{i,p,I}^{\bs{\bar s}}&:=(t_i-1)^p\widehat{\hspace{1mm}(t-1)_{i}^{\bs{\bar s}}}\zeta_I\frac{d}{dt_i},&\partial_{i,j,p,p',I}^{\bs{\bar s}}&:=(t_i-1)^{p}(t_j-1)^{p'}\hspace{-3mm}\widehat{\hspace{3mm}(t-1)_{i,j}^{\bs{\bar s}}}\zeta_I\frac{d}{dt_i}~~(i\neq j),
\end{array}\end{equation}}
we have 

{\begin{equation}\label{MMM}\begin{array}{rl}
\partial_{i,j,p}^{\bs{\bar s}}\in& (\bs{S}^+\D')^{0}\cap (\bs{S}^+\D')^{(p\ep_j-\ep_i+\sum_{q\neq j}s_q\ep_q)},\\
\partial_{i,p,I}^{\bs{\bar s}}\in &(\bs{S}^+\D')^{\d_I}\cap (\bs{S}^+\D')^{((p-1)\ep_i+\sum_{q\neq i}s_q\ep_q)},\\
\partial_{j,p,k,I}\in& (\bs{S}^+\D')^{\d_I-\d_k}\cap (\bs{S}^+\D')^{(p\ep_j)},\\
\partial_{i,j,p,p',I}^{\bs{\bar s}}\in& (\bs{S}^+\D')^{\d_I}\cap (\bs{S}^+\D')^{((p-1)\ep_i+p'\ep_j+\sum_{q\neq i, j}s_q\ep_q)}.
\end{array}
\end{equation}}

\begin{Pro}\label{power-2}
Recall $\fh$ and $\fh'$ from (\ref{h&h'}). Suppose that  $\bs{\theta}:\bs{S}^+\D'\longrightarrow \End(M)$ is a representation of $\bs{S}^+\D'$ in $M$. Assume the $\bs{S}^+\D'$-module $M$ has a weight space  decomposition $M=\Bigop{\mu\in\fh^*}{}M^\mu$ with respect to $\fh=\sspan_\bbbc\{h_k=\zeta_k{\partial_k}\mid 1\leq k\leq n\}$ with finite dimensional weight spaces.  
\begin{itemize}
\item[(i)] For $\lam\in (\fh')^*$, define 
\[M^{(\lam)}:=\{v\in M\mid \forall h'\in\fh' ~\exists p\in\bbbz^{>0}\hbox{~s.t.~}(\bs{\theta}(h')-\lam(h'){\rm id})^pv=0\}.\] Then, $M=\Bigop{\lam\in(\fh')^*}{}M^{(\lam)}$. Moreover, we have 
\[(\bs{S}^+\D')^{(\nu)}M^{(\lam)}\sub M^{(\lam+\nu)}\andd (\bs{S}^+\D')^{\upsilon}M^\mu\sub M^{\mu+\upsilon},\] 
for $\lam,\nu\in(\fh')^*$ and $\mu,\upsilon\in \fh^*,$
and that $M^\mu=\Bigop{\lam\in (\fh')^*}{}(M^\mu)^{(\lam)}$ with $(M^\mu)^{(\lam)}=M^\mu\cap M^{(\lam)}$. 
\item[(ii)] Suppose that weight multiplicities of $M=\op_{\mu\in\fh^*}M^\mu$ are bounded by a positive integer $d$, then,  for  $\ell:=(m+1)(2d+1)+1,$ we have  $(\bs{S}^+)^{\ell+n}\D' M=\{0\}$. 
\end{itemize}
\end{Pro}
\pf (i)  Let $\mu\in\fh^*.$ Then,  $h'_i=(t_i-1)\frac{d}{dt_i}$'s are $m+1$ commuting endomorphisms on the finite dimensional vector space $M^\mu.$ So, $M^\mu$ is decomposed into common generalized eigenspaces, that is,  $M^\mu=\Bigop{\lam\in(\fh')^*}{}(M^\mu)^{(\lam)}$ with 
\[(M^\mu)^{(\lam)}:=\{v\in M^\mu\mid \forall h'\in\fh' ~\exists p\in\bbbz^{>0}\hbox{~s.t~}(\bs{\theta}(h')-\lam(h'){\rm id})^pv=0\}.\] Since $M=\Bigop{\mu\in\fh^*}{}M^\mu,$ we get that $M=\Bigop{\lam\in(\fh')^*}{}M^{(\lam)}$ as we desired. For $\lam,\nu\in(\fh')^*,$ $x\in(\bs{S}^+\D')^{(\nu)}$ and  $h'\in\fh'$, we have 
\begin{align*}
{\small(\bs{\theta}(h')-(\lam+\nu)(h'){\rm id})\bs{\theta}(x)}
=&{\small\bs{\theta}(x)\bs{\theta}(h')-\lam(h')\bs{\theta}(x)=\bs{\theta}(x)(\bs{\theta}(h')-\lam(h'){\rm id}),}
\end{align*}  which in turn implies that $(\bs{S}^+\D')^{(\nu)}M^{(\lam)}\sub M^{(\lam+\nu)}.$ The other  assertions are clear facts in weight module theory.

(ii) It is enough to assume $\mu\in \fh^*,$ $\lam\in (\fh')^*$ as well as  $\bar{\bs{s}}:=(s_0,\ldots,s_m)\in(\bbbz^{\geq 0})^{m+1}$ with $s_0+\cdots+s_m>(m+1)(2d+1)$ and show that for  $v_0\in (M^\mu)^{(\lam)},$
\begin{itemize}
\item[(1)]  $(t_0-1)^{s_0}\cdots(t_m-1)^{s_m}\zeta_I\frac{d}{dt_i}v_0=0 \andd$
\item[(2)]  $(t_0-1)^{s_0}\cdots(t_m-1)^{s_m}\zeta_I{\partial_k}v_0=0,$
\end{itemize}
where  $0\leq i\leq m$,  $1\leq k\leq n$ and $I\sub\{1,\ldots,n\}.$

(1) To the contrary, assume $(t_0-1)^{s_0}\cdots(t_m-1)^{s_m}\zeta_I\frac{d}{dt_i}v_0\neq 0.$ Since $s_0+\cdots+s_m>(m+1)(2d+1)$,  there is $0\leq j\leq m$ such that  $s_j>2d+1.$  Recall (\ref{new-MM-1}).
If $i=j,$ for each $1\leq p\leq s_j$, we have 
\begin{align*}
(s_j+1-2p)(t_0-1)^{s_0}\cdots(t_m-1)^{s_m}\zeta_I\frac{d}{dt_i}v_0=&[\partial_{j,j,p}^{\bs{\bar{0}}},\partial^{\bs{\bar s}}_{j,s_j-p+1,I}]v_0\\
=&\partial_{j,j,p}^{\bs{\bar{0}}}(\partial^{\bs{\bar s}}_{j,s_j-p+1,I}v_0)-\partial^{\bs{\bar s}}_{j,s_j-p+1,I}(\partial_{j,j,p}^{\bs{\bar{0}}}v_0).
\end{align*}
In particular, if $2p\neq s_j+1,$ we have using  (\ref{MMM}) and part~(i) that either \[0\neq \partial_{j,j,p}^{\bs{\bar{0}}}v_0\in (M^{\mu})^{(\lam+p\ep_j-\ep_j)}\quad\hbox{or}\quad0\neq \partial^{\bs{\bar s}}_{j,s_j-p+1,I}v_0\in (M^{\mu+\d_I})^{(\lam+(s_j-p)\ep_j+\sum_{j\neq q}s_q\ep_q)}.\]  Since $s_j>2d+1,$ this implies that either $\dim(M^{\mu})\geq d+1$ or $\dim(M^{\mu+\d_I})\geq d+1$ which is a contradiction.
Also, if  $i\neq j,$ for $1\leq p\leq s_j$, we have 

\small{
\begin{align*}
0\neq (s_i+1)(t_0-1)^{s_0}\cdots(t_m-1)^{s_m}\zeta_I\frac{d}{dt_i}v_0
=&[\partial_{i,j,p}^{\bs{\bar{0}}},\partial^{\bs{\bar s}}_{i,j,s_i+1,s_j-p,I}]v_0\\
=&
\partial^{\bs{\bar{0}}}_{i,j,p}\partial^{\bs{\bar s}}_{i,j,s_i+1,s_j-p,I}v_0-\partial^{\bs{\bar s}}_{i,j,s_i+1,s_j-p,I}\partial^{\bs{\bar{0}}}_{i,j,p}v_0,
\end{align*}}
which implies that either  \[0\neq \partial_{i,j,s_i+1,s_j-p,I}^{\bs{\bar s}}v_0\in (M^{\mu+\d_I})^{(\lam+(s_j-p)\ep_j+\sum_{j\neq q}s_q\ep_q)}\quad\hbox{or}\quad0\neq \partial_{i,j,p}^{\bs{\bar{0}}}v_0\in (M^\mu)^{(\lam+p\ep_j-\ep_i)}.\] Again, we conclude that either $\dim(M^{\mu+\d_I})\geq d+1$ or $\dim(M^{\mu})\geq d+1$ which is a contradiction.
\smallskip

(2) To the contrary, assume $(t_0-1)^{s_0}\cdots(t_m-1)^{s_m}\zeta_I{\partial_k}v_0\neq 0.$  Since $s_0+\cdots+s_m>(m+1)(2d+1)$,  there is $0\leq j\leq m$ such that  $s_j>2d+1.$ For $1\leq p\leq s_j,$ we  have 
\begin{align*}
&\small{0\neq (s_j+1-p)(t_0-1)^{s_0}\cdots(t_m-1)^{s_m}\zeta_I{\partial_k}v_0=}\small{[(t_j-1)^{p}\widehat{(t-1)_j^{\bar s}}\frac{d}{dt_j},(t_j-1)^{s_j+1-p}\zeta_I{\partial_k}]v_0}\\
=&\small{[\partial^{\bs{\bar s}}_{j,j,p},\partial_{j,s_j+1-p,k,I}]v_0}=\small{\partial^{\bs{\bar s}}_{j,j,p}\hspace{-6mm}\underbrace{\partial_{j,s_j+1-p,k,I}v_0}_{\in (M^{\mu+\d_I-\d_k})^{(\lam+(s_j+1-p)\ep_j)}}\hspace{-6mm}-\hspace{6mm}\partial_{j,s_j+1-p,k,I}\hspace{-14mm}\underbrace{\partial^{\bs{\bar s}}_{j,j,p}v_0,}_{\in (M^{\mu})^{(\lam+(p-1)\ep_j+\sum_{j\neq q}s_q\ep_q)}}}
\end{align*}
which gives a contradiction as above.
\qed

\begin{Thm}\label{2-simple}
Suppose that $V$ is a simple $(\bs{S}\D'=)\bs{S}\D$-module such that 
\begin{itemize}
\item  $V$ has a weight space decomposition $V=\Bigop{\mu\in\fh^*}{}V^\mu$ with respect to $\fh=\sspan_\bbbc\{\zeta_i\partial_i\mid 1\leq i\leq n\}$,  
\item $\fh$-weight spaces are finite dimensional with bounded weight multiplicities.
\end{itemize}
Then,   $(\bs{S}^+)^2\D' V=\bs{S}^2\D V=\{0\}.$
\end{Thm}
\pf  Consider $V$ as an $\bs{S}^+\D'$-module. Then, by Proposition~\ref{power-2}, there is a positive integer $k$ such that 
\((\bs{S}^+)^k\D' V=\{0\}\) and so using Lemma~\ref{power-1}, we have 
\begin{equation}\label{zero-fi}
\bs{S}^{k+2}\D' V=\{0\}.
\end{equation}
Since $V$ is a simple module over $\bs{S}\D=\bs{S}\D'$, it is a simple module over $\bs{S}\D'/\bs{S}^{k+2}\D'\simeq \bs{S}^+\D'/(\bs{S}^+)^{k+2}\D'.$  Assume $\bar{ }:\bs{S}^+\D' \longrightarrow \bs{S}^+\D'/(\bs{S}^+)^{k+2}\D'$ is the canonical projection map and recall
\[D=\sum_{i=0}^m(t_i-1)\frac{d}{dt_i}+\sum_{k=1}^n\zeta_k{\partial_k}.\] Since $D$ is diagonalizable on $\bs{S}^+\D'$, $\bar D$ is diagonalizable on $\bs{S}^+\D'/(\bs{S}^+)^{k+2}\D'$.

Fix $\mu\in\fh^*$ in the support of $V=\Bigop{\lam\in\fh^*}{}V^\lam$ and  set  $\Omega:=V^\mu$. Assume $v\in \Omega$ is an eigenvector of $D$ corresponding to an  eigenvalue $\varrho.$ Since $V$ is a simple module over  $\bs{S}^+\D'/(\bs{S}^+)^{k+2}\D'$, its submodule generated by $v$ is the entire $V$. 
This together with the fact that 
 \[\bar D\cdot v=D\cdot v=\varrho v,\] implies that   $\bar D$ is diagonalizable on $V,$  say e.g., with the set of eigenvalues $P\sub \varrho+\bbbz^{\geq 0}$ and  $\varrho\in P.$ 
So, \[V=\Bigop{p\in P}{}V^p\quad  \hbox{with}\quad \{0\}\neq V^p=\{u\in V\mid D\cdot u=\bar D\cdot u=pu\}\quad(p\in P).\] This shows that  $D$ is diagonalizable on $V$ with the set of eigenvalues $P$. Since $\frac{(\bs{S}^+)^2\D'}{(\bs{S}^+)^{k+2}\D'}$ is an ideal of $\frac{\bs{S}^+\D'}{(\bs{S}^+)^{k+2}\D'}$ and $V$ is a simple module over  $\frac{\bs{S}^+\D'}{(\bs{S}^+)^{k+2}\D'}$, we get either 
\[(\bs{S}^+)^2\D'\cdot V=\frac{(\bs{S}^+)^2\D'}{(\bs{S}^+)^{k+2}\D'}\cdot V=\{0\}\quad\hbox{or}\quad (\bs{S}^+)^2\D'\cdot V=\frac{(\bs{S}^+)^2\D'}{(\bs{S}^+)^{k+2}\D'}\cdot V=V.\]
In the latter case, we get that $\varrho$ is an eigenvalue of $D$ on $(\bs{S}^+)^2\D'V$ which is impossible because  $P\sub\varrho+\bbbz^{\geq 0}$ and all eigenvalue of $D$ on $(\bs{S}^+)^2\D'$ are positive integers. So  we have $(\bs{S}^+)^2\D'\cdot V=\frac{(\bs{S}^+)^2\D'}{(\bs{S}^+)^{k+2}\D'}\cdot V=\{0\}$.
Finally, using Lemma~\ref{power-less} together with (\ref{zero-fi}), we get 
\[\bs{S}^2\D \cdot V=\bs{S}^2\D' \cdot V\sub (\bs{S}^+)^2\D' \cdot V +\bs{S}^{k+2}\D' \cdot V=\{0\}\] as we desired.\qed

\medskip

One  can easily get   the following lemma from (\ref{mods2}):
 \begin{lem}\label{gl}
We have ${\bs{S}}\D/\bs{S}^2\D\simeq \frak{gl}({m+1},n).$ In fact, \begin{align*}
&\theta:\bs{S}\D/\bs{S}^2\D\longrightarrow \frak{gl}(m+1,n)\\
&\left\{\begin{array}{rlrl}
(t_i-1)\bddtj+\bs{S}^2\D&\mapsto E_{i,j},
&\zeta_k\partial_l+\bs{S}^2\D&\mapsto E_{\bar k,\bar l},\\
(t_i-1)\partial_{k}+\bs{S}^2\D&\mapsto E_{i,\bar k},&\zeta_k\bddti+\bs{S}^2\D&\mapsto E_{\bar k,i},
\end{array}\right.
\end{align*}
is an isomorphism in which $0\leq i,j\leq m$ and $1\leq k,l\leq n.$
\end{lem}

\section{QP-Modules}
For a  Lie-Rinehart pair $(A,\fk)$,   the superspace  $\fg:=A\op \fk$ together with 
\begin{align*}
&(a\tinyop \d)\cdot(b\tinyop \sg):=a b\tinyop(a \sg+(-1)^{|b||x|}b \d) \quad(a,b\in A\andd \d,\sg\in \fk),
\end{align*}
is a supercommutative  associative superalgebra and it is a Lie superalgebra whose Lie bracket is defined as follows:
\begin{align*}
&\{a \tinyop \d,b \tinyop \sg\}:=(\d(b)-(-1)^{|a||\sg|}\sg(a)) \tinyop [\d,\sg] \quad(a,b\in A\andd \d,\sg\in \fk).
\end{align*}   
Moreover, $\fg$ together with the associative algebra homomorphism  $\pi:\fg\longrightarrow \fg$ mapping $a\tinyop \d$ to $a$,  is a {\it quasi-Poisson} superalgebra (or simply a QP-superalgebra) as defined in \cite{B}. We call  $\fg$  the {\it QP-superalgebra associated} to $(A,\fk).$

\begin{lem} \label{MM-40}
Suppose $({\dot A},{\dot \fk})$ is a Lie-Rinehart pair with corresponding QP-superalgebra $\fg={\dot A}\op{\dot \fk}$ and 
let $\pi:\fg\longrightarrow {\dot A}$  be the canonical projection map. 
Then, $(L({\dot A}),L(\fg))$ together with 
\[\begin{array}{ll}
(t_0^r\ot a)\cdot(t_0^s\ot b):=t_0^{r+s}\ot a\cdot b,&[t_0^r\ot x,t_0^s\ot y]:=t_0^{r+s}\ot (\{x,y\}{-r}x\cdot \pi(y){+s}\pi(x)\cdot y),\\\\
\begin{array}{c}
(t_0^s\ot a)(t_0^r\ot x):=t_0^{r+s}\ot a\cdot x,\\
(\hbox{\tiny{action of $L(\dot A)$ on $L(\fg)$}})
\end{array}
&\begin{array}{c}
\big(t_0^r\ot x\big)(t_0^s\ot b):=t_0^{r+s}\ot (s\pi(x)\cdot b+\{x,b\}),\\
(\hbox{\tiny{action of $L(\fg)$ on $L(\dot A)$}})
\end{array}
\end{array}\]
for $r,s\in\bbbz, a\in {\dot A}$ and $x\in \fg,$ is a Lie-Rinehart pair. In particular, if ${\dot A}:=\bbbc[{t}^{\pm1}_1,\ldots,{t}^{\pm1}_m]$ and ${\dot \fk}:={\rm Der}({\dot A}).$ The Lie superalgebra $L(\fg)$ is isomorphic to ${\rm Der}(L({\dot A}))$ in such a way that $t_0^ra(t_0\frac{d}{dt_0})\in {\rm Der}(L({\dot A}))$ {corresponds} to  $t_0^r\ot a\in L(\fg)$ and $t_0^r\d\in {\rm Der}(L({\dot A}))$ {corresponds} to  $t_0^r\ot \d\in L(\fg)$ for $a\in {\dot A}$ and $\d\in {\dot \fk}.$
\end{lem}

\begin{exa}[{Shen-Larsson Module}]\label{shen}
{\rm Suppose that $A:=\bbbc[t_0^{\pm1},\ldots,t_m^{\pm1}]\ot\Lam_n$
and $\fk:={\rm Der}(A).$ Let  $\Omega$ be a $\frak{gl}(m+1,n)$-module and  ${\bs{\bar \mu}}=(\mu(0),\mu(1),\ldots,\mu(m+n))\in\bbbc^{m+n+1}$ with 
\[\mu(m+k)=0\quad (1\leq k\leq n).\]
Set
\[\partial_\a:=\left\{
\begin{array}{ll}
\rmd_i& 0\leq \a=i\leq m,\\
\partial_{\a-m-1}& m+1\leq \a\leq m+n.
\end{array}
\right.\]
The superspace 
$V=V_A(\Omega,\bs{\bar \mu}):=A\ot \Omega$ together with 
\begin{align*}
f\partial_\a* (g\ot v)=&\sum_{\gamma=0}^{m+n}(-1)^{|\partial_\gamma|+|fg||\partial_\gamma|+|g||\partial_\a|}\partial_\gamma(f)g\ot E_{\gamma,\a}\omega+f(\partial_\a(g)+\mu(\a) g)\ot \omega,\\
f\cdot(g\ot \omega)=&fg\ot\omega,
\end{align*}
is an $A\fk$-module. This module is called a Shen-Larsson module. If $\Omega$ is simple, then $V_A(\Omega,\bs{\bar \mu})$ is simple; see \cite[Lem.'s 2.2, 3.3 $\&$ 3.5(1)]{XL1}
}
\end{exa}

\begin{deft}\label{qp-mod}
{\rm Suppose that $(A,\fk)$ is a  Lie-Rinehart pair  with associated QP-superalgebra $(\fg,~\cdot~,\{\cdot,\cdot\})$.  We say a superspace $\Omega$ is a {\it quasi-Poisson  module} (or simply a {\it QP-module}) over the Lie-Rinehart pair $(A,\fk)$  if there are 
\begin{align*}
\varphi:&A\longrightarrow \End(\Omega),\\
\psi:&\fg\longrightarrow \End(\Omega),\\
\hat \varphi:&\fg\longrightarrow \End(\Omega),
\end{align*}
such that 
\begin{itemize}
\item[(1)] $\varphi$ is an associative superalgebra homomorphism, i.e., $\Omega$ is an $A$-module,
\item[(2)] $\psi$ is a Lie superalgebra homomorphism, i.e., $\Omega$ is a $(\fg,\{\cdot,\cdot\})$-module,
\item[(3)]  $[\hat \varphi_x ,\hat  \varphi_y]=\hat  \varphi_{x\cdot \pi(y)}-\hat  \varphi_{\pi(x)\cdot y},$ 
\item[(4)] $\hat\varphi_{a\cdot x}=\varphi_a\hat\varphi_x,$
\item[(5)] $[\hat\varphi_x,\varphi_a]=0,$ 
\item[(6)] $[\hat  \varphi_x ,\psi_y]=\hat  \varphi_{\{x,y\}}-(-1)^{|x||y|}\varphi_{_{\pi(y)}}\psi_{_x}+\psi_{_{x\cdot \pi(y)}},$
\item[(7)] $[\psi_{_x},\varphi_{a}]=\varphi_{x(a)},$
\end{itemize}
for $a\in A$ and $x,y\in \fg;$ here $x(a)$ means $\{x,a\}.$ 
A linear homomorphism $\Theta$ from a QP-module $(\Omega,\varphi,\hat\varphi,\psi)$ to  a QP-module $(\dot \Omega,\dot \varphi,\hat{\dot\varphi},\dot \psi)$ is  called a homomorphism if 
\begin{align*}
\Theta(\psi_x(\omega))=\dot \psi_x(\Theta(\omega)), \Theta(\varphi_a(\omega))=\dot\varphi_a(\Theta(\omega))\andd \Theta(\hat\varphi_x(\omega))=\hat{\dot\varphi}_x(\Theta(\omega))
\end{align*}
for $x\in\fg,$ $a\in A$ and $\omega\in \Omega.$
The notions {\it submodule}, {\it simplicity} and {\it isomorphism},  are defined in the usual manner.}
\end{deft}
\begin{exa}
{\rm Suppose that $(A,\fk)$ is a  Lie-Rinehart pair  with associated QP-superalgebra $(\fg,~\cdot~,\{\cdot,\cdot\})$ with the canonical projection map $\pi:\fg\longrightarrow A.$ Define \begin{align*}
\varphi:&A\longrightarrow \End(\fg),
&&\varphi_a(y)=a\cdot y,\\
\psi:&\fg\longrightarrow \End(\fg),
&&\psi_x(y)=\{x,y\},\\
\hat \varphi:&\fg\longrightarrow \End(\fg),
&&\hat\varphi_x(y)=x\cdot\pi(y)=(-1)^{|x||y|}\pi(y)\cdot x,
\end{align*}
for $a\in A$ and  $x,y\in \fg$. Then, $\fg$ together with $(\varphi,\hat\varphi,\psi)$ is a QP-module over $(A,\fk).$}
\end{exa}
\begin{lem}\label{MM-10}
Suppose that $(M,\varphi,\hat \varphi,\psi)$ is a QP-module over a Lie-Rinehart pair $(A,\fk)$ with corresponding QP-superalgebra  $\fg.$ Then 
$L(M)=\bbbc[t_0^{\pm1}]\ot M$ together with 
\begin{align*}
(t_0^r\ot x)*(t_0^s\ot \omega):=&t_0^{r+s}\ot (\psi_{_x}(\omega){-r}\hat\varphi_x(\omega){+s}\varphi_{\pi(x)}(\omega)),\\
(t_0^r\ot a)\cdot (t_0^s\ot \omega):=&t_0^{r+s}\ot \varphi_a(\omega),
\end{align*}
is an  $L(A)L(\fg)$-module.
 Moreover, if the $L(A)L(\fg)$-module $L(M)$ is simple, then, $M$ is a simple QP-module.
\end{lem}
\pf 
It is not hard to see that $L(M)$ is an  $L(A)L(\fg)$-module. For the last assertion, we first note that if $v,w\in M$ such that $v$ is an element of QP-submodule generated by $w$, then 
${(t_0^p\ot a)}\cdot (t_0^s\ot v)$ and ${(t_0^r\ot x)*(t_0^s\ot v)}$ are respectively of the form $t_0^{p+s}\ot u_1$ and  $t_0^{r+s}\ot u_2$ for some  elements $u_1$ and $u_2$  of the QP-submodule  generated by $w$.

Now, suppose that $L(M)$ is simple and assume $v\in M.$  Since $v$ is an element of the QP-submodule generated by $v$,  each element of $L(A)L(\fg)$-submodule generated by 
$v$ (in fact each element of  $L(M)$ due to the simplicity of  $L(M)$) is a linear combinations of elements of the form $t_0^r\ot u$ where $u$ is an element of QP-submodule generated by $v.$ In particular, if $\omega\in M\sub L(M),$ then $\omega=\sum_{i=1}^p t_0^{r_i}\ot u_i$ in which $u_i$ is an element of QP-submodule generated by $v.$ This implies that $\sum_{i=1,r_i\neq 0}^p t_0^{r_i}\ot u_i=0$ and $\omega$ is an element of QP-submodule generated by $v.$
\qed

\begin{exa}{\rm
Set   ${\dot A}:=\bbbc[t_1^{\pm1},\ldots,t_m^{\pm1}]\ot \Lam(\zeta_1,\ldots,\zeta_n)$ and ${\dot \fk}:={\rm Der}(\dot A)=\Bigop{i=1}{m}\dot A\bddti\op\Bigop{k=1}{n}\dot A\partial_{k}.$  Suppose that $\Omega$ is a  $\frak{gl}(m+1,n)$-module and ${\bs{\bar \mu}}=(\mu(0),\mu(1),\ldots,\mu(m+n))\in\bbbc^{m+n+1}$ with 
\[\mu(m+k)=0\quad (1\leq k\leq n).\]
For $0\leq \a\leq m+n+1,$ set 
\[\partial_\a:=\left\{
\begin{array}{ll}
1&\a=0,\\
\rmd_i&1\leq \a=i\leq m,\\
\partial_k& m+1\leq \a=\bar k=m+k\leq m+n.
\end{array}
\right.
\]Next set
\[M_{\dot A}(\Omega,\bs{\bar\mu}):=\dot A\ot_\bbbc \Omega\] and  define  

\begin{align*}
\varphi_a(b\ot v):=&a b\ot v,\\
\hat \varphi_{a\partial_\a}(b\ot v):=&-(-1)^{|b||\partial_\a|}ab\ot E_{0,\a}v,\\
\psi_{{a\partial_\a}}(b\ot v):=&a(\partial_\a(b)+\mu(\a) b)\ot v+\sum_{\b=1}^{m+n}{(-1)^{|ab||\partial_\b|+|\partial_\b|+|b||\partial_\a|}}\partial_\b(a)b\ot E_{_{\b,\a}}v,
\end{align*}
for  $a,b\in {\dot A}$ and $v\in \Omega.$ Then, $M_{\dot A}(\Omega,\bs{\bar \mu})$ is a  QP-module  over the Lie-Rinehart pair $({\dot A},{\dot \fk})$ and we have $V_A(\Omega,\bs{\bar \mu})\simeq L(M_{\dot A}(\Omega,\bs{\bar \mu}));$ see Example~\ref{shen} and Lemma~\ref{MM-10}. We call $M_{\dot A}(\Omega,\bs{\bar \mu})$ a Shen-Larsson QP-module. Also, if $\Omega$ is simple, then  $V_A(\Omega,\bs{\bar \mu})=L(M_{\dot A}(\Omega,\bs{\bar \mu}))$ is a  simple  $A\fk$-module and so
$M_{\dot A}(\Omega,\bs{\bar \mu})$ is simple using Lemma~\ref{MM-10}.
}
\end{exa}

{From now until end of this section,} we assume  \[{\dot A}:=\bbbc[t_1^{\pm1},\ldots,t_m^{\pm1}]\ot \Lam(\zeta_1,\ldots,\zeta_n)\andd {\dot \fk}:={\rm Der}({\dot A})=\Bigop{i=1}{m}{\dot A}\bddti\op\Bigop{k=1}{n}{\dot A}\partial_{k},\] and that   \[A:=\bbbc[t_0^{\pm1},\ldots,t_m^{\pm1}]\ot \Lam(\zeta_1,\ldots,\zeta_n)\andd \fk:={\rm Der}(A)=\Bigop{i=0}{m}A\bddti\op\Bigop{k=1}{n}A\partial_{k}.\]

The  subspace ${A}\smpro   \fk$ of ${A}\smpro   U(\fk)$ is a subalgebra of the Lie superalgebra of  the associative superalgebra ${A}\smpro   U(\fk).$ It is a $\bbbz^{m+1}$-graded Lie superalgebra with  
\[({A}\smpro   {\fk})^{(r_0,\ldots,r_m)}=\sum t^{\bs{\bar\b}}\zeta_I\smpro   t^{\bs{\bar\gamma}}\zeta_{J}\d\]
where the summation runs over all $\d\in\{\rmd_i,{\partial_k} \mid {0\leq i\leq m, 1\leq k\leq n}\},$ $I,J\sub\{1,\ldots,n\}$ and   $\bs{\bar\b},\bs{\bar\gamma}\in\bbbz^{m+1}$  with  \((r_0,\ldots,r_m)=\bs{\bar\b}+\bs{\bar\gamma}.\)

\begin{lem}[{\cite[\S~3]{XL1}}]\label{XL1}For  subsets\footnote{$I$ or $J$ is an  empty set if $p=0$ or $q=0$ respectively.} $I=\{i_1<\cdots<i_p\}$ and $J=\{j_1<\cdots<j_q\}$ of $\{1,\ldots,n\}$, set 
\[\tau(I,J):=\sum_{s=1}^q|\{r\mid j_s<i_r\}|\]
in which by $|\cdot|,$ we mean the cardinal number of a set. Also, we set
\begin{equation}\label{XX}
X_{\bs{\bar{r}},J,\partial}:=\left\{
\begin{array}{ll}
\displaystyle{\sum_{I\sub J}(-1)^{|I|+\tau(I,J\setminus I)}t^{-{\bs{\bar{r}}}}\zeta_{_I}\smpro   t^{\bs{\bar{r}}} \zeta_{_{J\setminus I}}\partial} &\hbox{ if $J\neq \emptyset,$}\\
(t^{-{\bs{\bar{r}}}}\smpro   t^{\bs{\bar{r}}} \partial)-\partial & \hbox{if $J=\emptyset,$}
\end{array}\right.
\end{equation}
for $\bs{\bar r}\in\bbbz^{m+1}$ and  \(\partial\in {\{\rmd_i,\partial_k\mid 0\leq i\leq m, 1\leq k\leq n}\}\).
\begin{itemize}
\item[(i)] Assume $\T$ is a subspace  of $A\smpro   \fk $ spanned by
$X_{\bs{\bar{r}},I,\partial},X_{\bs{\bar{s}},\emptyset,\partial}$ where $ \bs{\bar r}=(r_0,\ldots,r_m),\bs{\bar s}=(s_0,\ldots,s_m)$ run over $\bbbz^{m+1}$ with $\bs{\bar s}\neq 0$, $\emptyset\neq I\subset \{1,\ldots,n\}$ and $\partial\in \{\rmd_i,\partial_k\mid 0\leq i\leq m, 1\leq k\leq n\}.$ 
Then, recalling (\ref{deltadelta'}), we have  
\[\T=\{x\in (A\smpro \fk)^0\mid x\smpro  \d=\d\smpro  x,x\smpro  a=a\smpro  x~(a\in A,\d\in \D)\}.\] In particular, $\T$ is a subalgebra of the Lie superalgebra $(A\smpro \fk)^0$.
\item[(ii)] The Lie superalgebra $\T$ is isomorphic to $\bs{S}\D$ (see  \S~\ref{section1}) via the isomorphism 
\[\begin{array}{llll}
\Psi:X_{\bs{\bar{r}},\emptyset,\partial}\mapsto (t^{\bs{\bar{r}}}-1)\partial, &X_{\bs{\bar{r}},I,\partial}\mapsto t^{\bs{\bar{r}}}\zeta_I\partial,
\end{array}\] 
in particular, 
\[\Psi(t_i\smpro  \frac{d}{dt_i}-t_i\frac{d}{dt_i})={-(t_i-1)\frac{d}{dt_i}}\andd \Psi(\zeta_k\smpro {\partial_k}-\zeta_k{\partial_k})={-\zeta_k{\partial_k}}.\]
\end{itemize}
\end{lem}

\begin{lem} \label{equalities} Assume $(M,\varphi,\psi,\hat\varphi)$ is a QP-module over $({\dot A},{\dot \fk})$. For $a,b\in {\dot A}$, $I\sub\{1,\ldots,n\}$, $\bs{\bar r}\in \bbbz^m$ and $\partial\in\fg,$ we have the following:
\begin{itemize}
\item[(1)]  $\psi_{ab}=-\varphi_{ab}\psi_1+\varphi_a\psi_{ b}+(-1)^{|a||b|}\varphi_{b}\psi_a.$
\item[(2)] $\hat\varphi_{\partial(a)}=\hat\varphi_\partial\psi_a-(-1)^{|a||\partial|}\psi_a\hat\varphi_\partial+(-1)^{|a||\partial|}\varphi_a\psi_\partial-(-1)^{|a||\partial|}\psi_{a\partial}.$
\item[(3)] $\psi_{t^{\bs{\bar r}}a}=\varphi_{t^{\bs{\bar r}}}\psi_a+\sum_{j=1}^m\varphi_{a\rmd_j(t^{\bs{\bar r}})}(\varphi_{t_j^{-1}}\psi_{t_j}-\psi_1).$
\item[(4)] $\psi_{\zeta_I}=\varphi_{\zeta_I}\psi_1+(-1)^{|I|}\sum_{p=1}^n\varphi_{\partial_p(\zeta_I)}(-\psi_{\zeta_{p}}+\varphi_{\zeta_{p}}\psi_1).$
\item[(5)] $\psi_{a\partial}=(-1)^{|a|+1}\sum_{p=1}^n\varphi_{\partial_p(a)}(\psi_{\zeta_{p}\partial}-\varphi_{\zeta_{p}}\psi_{\partial})
+\sum_{{j=1}}^m\varphi_{\rmd_j(a)}
(\varphi_{t_{j}^{-1}}\psi_{t_{j} \partial}-\psi_{ {\partial}})+\varphi_{a}\psi_{\partial}.$
\end{itemize}
\end{lem}
\pf (1) We have 
{\footnotesize\begin{align*}
&\psi_{ab}=[\hat  \varphi_a ,\psi_b]+(-1)^{|a||b|}\varphi_{b}\psi_a=\hat  \varphi_a \psi_b-(-1)^{|a||b|}\psi_b\hat  \varphi_a+(-1)^{|a||b|}\varphi_{b}\psi_a\\
=&\varphi_a\hat  \varphi_1 \psi_b-(-1)^{|a||b|}\psi_b\hat  \varphi_a+(-1)^{|a||b|}\varphi_{b}\psi_a
=\varphi_a\psi_b\hat  \varphi_1-\varphi_a\varphi_{b}\psi_1+\varphi_a\psi_{ b} -(-1)^{|a||b|}\psi_b\hat  \varphi_a+(-1)^{|a||b|}\varphi_{b}\psi_a\\
=&(-1)^{|a||b|}\psi_b\varphi_a\hat  \varphi_1-\varphi_a\varphi_{b}\psi_1+\varphi_a\psi_{ b} -(-1)^{|a||b|}\psi_b \varphi_a\hat  \varphi_1+(-1)^{|a||b|}\varphi_{b}\psi_a
=-\varphi_a\varphi_{b}\psi_1+\varphi_a\psi_{ b}+(-1)^{|a||b|}\varphi_{b}\psi_a.
\end{align*}}

(2) is just  condition (6) of  Definition~\ref{qp-mod}.

(3) For $1\leq i\leq m,$ we have 
\(\psi_{t_i}=\varphi_{t_i}\psi_1+\varphi_{t_i}(\varphi_{t_i^{-1}}\psi_{t_i}-\psi_1)\). This together with part~(1) 
implies that 
\begin{align*}
\psi_{t_i^{-1}}=\varphi_{t_i^{-1}}(\varphi_{t_i}\psi_{t_i^{-1}})
=\varphi_{t_i^{-1}}(\psi_1+\psi_1-\varphi_{t_i^{-1}}\psi_{t_i})
=\varphi_{t_i^{-1}}\psi_1-\varphi_{t_i^{-1}}(\varphi_{t_i^{-1}}\psi_{t_i}-\psi_1).
\end{align*}
Therefore, for  $a\in A$ and $\sg=\pm1,$ we have 

\begin{align*}
\psi_{t_i^\sg a}=&-\varphi_{t_i^\sg a}\psi_1+\varphi_{t_i^\sg}\psi_a+\varphi_a\psi_{t_i^\sg}=-\varphi_{t_i^\sg a}\psi_1+\varphi_{t_i^\sg}\psi_a+\varphi_a(\varphi_{t_i^\sg}\psi_1+\sg\varphi_{t_i^\sg}(\varphi_{t_i^{-1}}\psi_{t_i}-\psi_1))\\
=&\varphi_{t_i^\sg}\psi_a+\sg\varphi_{at_i^\sg}(\varphi_{t_i^{-1}}\psi_{t_i}-\psi_1)=\varphi_{t_i^\sg}\psi_a+\sum_{j=1}^m\varphi_{a\rmd_j(t_i^\sg)}(\varphi_{t_j^{-1}}\psi_{t_j}-\psi_1).
\end{align*}
Now, for $\bs{\bar{s}}=(s_1,\ldots,s_m),$ we use induction on $|\bs{\bar{s}}|:=|s_1|+\cdots+|s_m|$ to get the result: We have already seen the result for $|\bs{\bar{s}}|=1.$ We then assume $|\bs{\bar{s}}|>1$ and that the result holds for all $\bs{\bar r}$ with $|\bs{\bar r}|<|\bs{\bar s}|.$ We have $t^{\bs{\bar s}}=t_i^\sg t^{\bs{\bar r}}$ for some $1\leq i\leq m,$ $\sg=\pm1$ and $\bs{\bar{r}}=(r_1,\ldots,r_m)$ with $|\bs{\bar r}|=|\bs{\bar s}|-1.$ For $a\in A$, we have 
\begin{align*}
\psi_{t^{\bs{\bar s}}a}=&\psi_{t_i^\sg t^{\bs{\bar r}}a}=\varphi_{t_i^\sg}\psi_{{t^{\bs{\bar r}}a}}+\sum_{j=1}^m\varphi_{{t^{\bs{\bar r}}a}\rmd_j(t_i^\sg)}(\varphi_{t_j^{-1}}\psi_{t_j}-\psi_1)\\
=&\varphi_{t_i^\sg}\varphi_{t^{\bs{\bar r}}}\psi_a+\sum_{j=1}^m\varphi_{t_i^\sg}\varphi_{a\rmd_j(t^{\bs{\bar r}})}(\varphi_{t_j^{-1}}\psi_{t_j}-\psi_1)+\sum_{j=1}^m\varphi_{{t^{\bs{\bar r}}a}\rmd_j(t_i^\sg)}(\varphi_{t_j^{-1}}\psi_{t_j}-\psi_1)\\
=&\varphi_{t_i^\sg}\varphi_{t^{\bs{\bar r}}}\psi_a+\sum_{j=1}^m\varphi_{a\rmd_j(t_i^\sg t^{\bs{\bar r}})}(\varphi_{t_j^{-1}}\psi_{t_j}-\psi_1)=\varphi_{t^{\bs{\bar s}}}\psi_a+\sum_{j=1}^m\varphi_{a\rmd_j(t^{\bs{\bar s}})}(\varphi_{t_j^{-1}}\psi_{t_j}-\psi_1).
\end{align*}

(4) follows easily from (1) and  an induction process.

(5) We have 
{\footnotesize\begin{align*}
&\psi_{a\partial}=  (-1)^{|a||\partial|}[\hat\varphi_\partial,\psi_a]-(-1)^{|a||\partial|}\hat\varphi_{\partial(a)}+\varphi_a\psi_\partial\\
=&(-1)^{|a||\partial|}[\hat\varphi_\partial,
\varphi_{a}\psi_1+(-1)^{|a|}\sum_{p=1}^n\varphi_{\partial_p(a)}(-\psi_{\zeta_{p}}+\varphi_{\zeta_{p}}\psi_1)+\sum_{j=1}^m\varphi_{\rmd_j(a)}(\varphi_{t_j^{-1}}\psi_{t_j}-\psi_1)]
-(-1)^{|a||\partial|}\hat\varphi_{\partial(a)}+\varphi_a\psi_\partial\\
=&-\sum_{p=1}^n(-1)^{|a|}(-1)^{|a||\partial|}(-1)^{|\partial||\partial_p(a)|}\varphi_{\partial_p(a)}[\hat\varphi_\partial,\psi_{\zeta_{p}}]+\sum_{j=1}^m\varphi_{\rmd_j(a)}\varphi_{t_j^{-1}}[\hat\varphi_\partial,\psi_{t_j}]
-(-1)^{|a||\partial|}\hat\varphi_{\partial(a)}+\varphi_a\psi_\partial\\
=&-\sum_{p=1}^n(-1)^{|a|+|\partial|}\varphi_{\partial_p(a)}(\hat\varphi_{\partial(\zeta_p)}-(-1)^{|\partial|}\varphi_{\zeta_p}\psi_\partial+(-1)^{|\partial|}\psi_{\zeta_p\partial})+\sum_{j=1}^m\varphi_{\rmd_j(a)}\varphi_{t_j^{-1}}(\hat\varphi_{\partial(t_j)}-\varphi_{t_j}\psi_\partial+\psi_{t_j\partial})\\
-&(-1)^{|a||\partial|}\hat\varphi_{\partial(a)}+\varphi_a\psi_\partial\\
=&\sum_{p=1}^n(-1)^{|a|}\varphi_{\partial_p(a)}(\varphi_{\zeta_p}\psi_\partial-\psi_{\zeta_p\partial})+
\sum_{j=1}^m\varphi_{\rmd_j(a)}\varphi_{t_j^{-1}}(-\varphi_{t_j}\psi_\partial+\psi_{t_j\partial})\\
-&\sum_{p=1}^n(-1)^{|a|+|\partial|}\varphi_{\partial_p(a)}\hat\varphi_{\partial(\zeta_p)}+\sum_{j=1}^m\varphi_{\rmd_j(a)}\varphi_{t_j^{-1}}\hat\varphi_{\partial(t_j)}
-(-1)^{|a||\partial|}\hat\varphi_{\partial(a)}+\varphi_a\psi_\partial\\
=&\sum_{p=1}^n(-1)^{|a|}\varphi_{\partial_p(a)}(\varphi_{\zeta_p}\psi_\partial-\psi_{\zeta_p\partial})+
\sum_{j=1}^m\varphi_{\rmd_j(a)}\varphi_{t_j^{-1}}(-\varphi_{t_j}\psi_\partial+\psi_{t_j\partial})+\varphi_a\psi_\partial.
\end{align*}}
This completes the proof.
\qed

{From now on, we fix  a simple QP-module $(M,\varphi,\hat\varphi,\psi)$ over $({\dot A},{\dot \fk})$}. By abuse of  notations, we 
set
\begin{equation}\label{conv-zero}
\hat\varphi_{\rmd_0}:=\hat\varphi_1\andd \psi_{{a}\rmd_0}:=\psi_a\quad(a\in \dot A).
\end{equation}
For $\bs{\bar r}=(r_0,\ldots,r_m)\in\bbbz^{m+1}$, set 
\begin{equation}\label{prim}\bs{{\bar r}'}=(r_1,\ldots,r_m)\in\bbbz^{m}.
\end{equation}
 Using  Lemmas~\ref{MM-40} \& \ref{MM-10}, $L(M)=\bbbc[t_0^{\pm1}]\ot M$ is an $A\fk$-module (equivalently an $(A\smpro U(\fk))$-module) with 
\[(t^{\bs{\bar{r}}}\zeta_I\smpro  t^{\bs{\bar{s}}}\zeta_J\partial)* (t_0^{k}\ot u):=t_0^{r_0+s_0+k}\ot \varphi_{t^{\bs{\bar r}'}\zeta_I}(\psi_{ t^{\bs{\bar{s}}'}\zeta_J\partial}(u)-s_0\varphi_{ t^{\bs{\bar{s}}'}\zeta_J}\hat\varphi_\partial(u)+k\d_{\rmd_0,\partial}\varphi_{ t^{\bs{\bar{s}}'}\zeta_J}(u))\] in which $u\in M$ and $\partial\in\{\rmd_i,\partial_k\mid 0\leq i\leq m, 1\leq k\leq n \}.$  We note that  
\begin{equation}\label{MM-T}
(t^{-\bs{\bar{r}}}\zeta_I\smpro  t^{\bs{\bar{r}}}\zeta_J\partial)* u= \varphi_{t^{-\bs{\bar r}'}\zeta_I}(\psi_{ t^{\bs{\bar{r}}'}\zeta_J\partial}(u)-r_0\varphi_{ t^{\bs{\bar{r}}'}\zeta_J}\hat\varphi_\partial(u)).
\end{equation}
In particular, for $\bs{\bar r}\in\bbbz^{m+1}$ and  \(\partial\in {\{\rmd_i,\partial_k\mid 0\leq i\leq m, 1\leq k\leq n}\}\),  recalling (\ref{XX}), we have 
\begin{equation}\label{MM-30}
\begin{array}{rl}
X_{\bs{\bar r},\emptyset,\partial}*u=& \varphi_{t^{-\bs{\bar r}'}}(\psi_{ t^{\bs{\bar{r}}'}\partial}(u)-r_0\varphi_{ t^{\bs{\bar{r}}'}}\hat\varphi_\partial(u)) -\psi_{\partial}(u),\\
X_{\bs{\bar r},J,\partial}*u=&
\displaystyle{\sum_{J'\sub J}(-1)^{|J'|+\tau(J',J\setminus J')}\varphi_{t^{-\bs{\bar r}'}\zeta_{J'}}(\psi_{ t^{\bs{\bar{r}}'}\zeta_{J\setminus J'}\partial}u-r_0
\varphi_{ t^{\bs{\bar{r}}'}\zeta_{J\setminus J'}}\hat\varphi_\partial(u))}\andd \\
X_{(1,0,\ldots,0),\emptyset,\partial}*u=&-\hat\varphi_\partial(u) \quad (u\in M).
\end{array}
\end{equation}
\begin{lem}\label{MM-Omega}
Assume $\{0\}\neq N\sub M$ is invariant under $\psi_{\partial_k}$ $(1\leq k\leq n)$.
Then, 
 \[\Omega(N):=\{u\in N\mid \psi_{\partial_k}u=0~~(1\leq k\leq n)\}\] is nonzero.  
\end{lem}
\pf 
We note that for $w\in M$ and  $1\leq k\leq n,$ we have  
\[\psi_{\partial_{k}}\psi_{\partial_{k}} w=\frac{1}{2}[\psi_{\partial_{k}},\psi_{\partial_{k}}]w=\frac{1}{2}\psi_{[\partial_{k},\partial_{k}]}w=0.\]
Fix $0\neq u\in N$ and consider the set 
\[\{\psi_{\partial_{j_1}}\cdots\psi_{\partial_{j_s}}u\neq 0\mid 0\leq s\leq n,j_1<\ldots<j_s\}\] in which for $s=0,$ we define $\psi_{\partial_{j_1}}\cdots\psi_{\partial_{j_s}}u:=u.$
Set
\[p_0:={\rm max}\{0\leq p\leq n\mid \exists 1\leq j_{1,p}<\cdots<j_{s,p}\leq n\hbox{ s.t. }u_p:=\psi_{\partial_{j_{1,p}}}\cdots\psi_{\partial_{j_{s,p}}}u\not=0\}.\] For $u_0:=u_{p_0}$, we have $u_0\in N$ and that  $\psi_{\partial_k} u_0=0$ for all $1\leq k\leq n.$
\qed

\begin{deft}
{\rm  The QP-module $(M,\varphi,\hat\varphi,\psi)$ is called {\it cuspidal} if $(M,\psi)$ has a weight space decomposition with respect to \[\ft:=\sspan_\bbbc\{1,\rmd_i,\zeta_k\partial_k\mid 1\leq i\leq m, 1\leq k\leq n\}\] with finite dimensional weight spaces and bounded weight multiplicities. }
\end{deft}
\begin{Thm} \label{main}
Assume the QP-module $(M,\varphi,\hat\varphi,\psi)$ is simple and cuspidal with $M=\Bigop{\lam\in\ft^*}{}M^\lam$. Set $\ft':=\sspan_\bbbc\{1,\rmd_i\mid 1\leq i\leq m\}$ and for $\mu\in(\ft')^*$, define 
\[M_{(\mu)}:=\bigoplus_{\substack{\lam\in\ft^*\\
\lam\mid_{\ft'}=\mu}}M^\lam.\]
Fix $\mu_0\in (\ft')^*$ and set
\[N:=M_{(\mu_0)}\andd \Omega:=\Omega(N);\] see Lemma~\ref{MM-Omega}.
We have the following:
\begin{itemize}
\item[(i)] The subset $\Omega$ is invariant under  $\varphi_{t_{j}^{-1}}\psi_{t_{j} \partial}-\psi_{ {\partial}}$,  $\varphi_{\zeta_{p}}\psi_{\partial}-\psi_{\zeta_{p}\partial}$ and $\hat\varphi_\partial$ for  $1\leq j\leq m,$ $1\leq p\leq n$ and $\partial\in\{1,\partial_k,\rmd_i\mid 1\leq k\leq n, 1\leq i\leq m\}$.
\item[(ii)]  {$\Omega$  is a $\T$-submodule of the $(A\smpro\fk)$-module $L(M)$ (see Lemma~\ref{XL1}(i)) having weight space decomposition with respect to \(\{\zeta_k\partial_k\mid 1\leq k\leq n\}\) with finite dimensional weight spaces and bounded weight multiplicities}.
\item[(iii)] If $U$ is a $\T$-submodule of $\Omega$, then \[\sspan_\bbbc\{\varphi_{\zeta_I}u\mid I\sub\{1,\ldots,n\},u\in U\}=N\andd \sspan_\bbbc\{\varphi_{a}u\mid a\in {\dot A},u\in U\}=M.\]
\item[(iv)] $\Omega$ is a simple $\T$-module. In particular, $\Omega$ is a simple admissible $\frak{gl}(m+1,n)$-module via the representation $\Phi$ defined by 
\[\begin{array}{lll}
\Phi(E_{0,0}):={-\hat\varphi_{_1}}&\Phi(E_{0,i}):={-\hat\varphi_{\rmd_i}}&\Phi(E_{0,\bar k}):=-\hat\varphi_{{\ddk}}\\
\Phi(E_{i,0}):={\varphi_{{t_i}^{-1}}\psi_{t_i}-\psi_1}&\Phi(E_{i,j}):={\varphi_{{t_i}^{-1}}\psi_{t_i\rmd_j}-}\psi_{\rmd_j}&\Phi(E_{i,\bar k}):={\varphi_{{t_i}^{-1}}\psi_{t_i\ddk}-}\psi_{\ddk}\\
\Phi(E_{\bar k,0}):={\psi_{\zeta_k}-\varphi_{\zeta_k}}\psi_1&{\Phi(E_{\bar k,i}):=\psi_{\zeta_k\rmd_j}-\varphi_{\zeta_k}\psi_{\rmd_j}}
&\Phi(E_{\bar k,\bar l}):=\psi_{\zeta_k\bddl}-\varphi_{\zeta_k}\psi_{\bddl}.
\end{array}
\] 
\item[(v)] {The simple  cuspidal QP-module $(M,\varphi,\hat\varphi,\psi)$ is a Shen-Larsson QP-module.}
\end{itemize}
\end{Thm}
\pf  (i) follows from the properties of a QP-module.

(ii) Using  (\ref{MM-T}), (\ref{MM-30}) and  Lemma~\ref{XL1}(i), we get that $\Omega$ is a $\T$-submodule of $L(M).$
The last assertion follows from the fact that $\Omega\sub N$ is invariant under  $\{\psi_{\zeta_k\partial_k}\mid 1\leq k\leq n\}$. 

(iii) Set $W:=\sspan_\bbbc\{\varphi_{a}u\mid a\in {\dot A},u\in U\}.$ We shall show that $W$ is a submodule of the QP-module $(M,\varphi,\hat\varphi,\psi).$ Since $U$ is invariant under $\T,$ by (\ref{MM-T}), it is invariant under 
\[\varphi_{\zeta_{p}}\psi_{\partial}-\psi_{\zeta_{p}\partial}\andd \varphi_{t_{j}^{-1}}\psi_{t_{j} \partial}-\psi_{ {\partial}}\]
for $1\leq p\leq n$, $1\leq j\leq m$ and $\partial\in\{1,\rmd_i,\partial_k\mid 1\leq i\leq m,1\leq k\leq n\}.$ So, 
for $u\in U$, $b\in {\dot A}$ and $\partial\in\{1,\rmd_i,\partial_k\mid 1\leq i\leq m,1\leq k\leq n\},$ using Lemma~\ref{equalities}, we have
\begin{align*}
\psi_{b\partial}\varphi_{a}u=&\varphi_{b\partial(a)}u+(-1)^{|a||b\partial|}\varphi_{a}\psi_{b\partial}u\\
=&\varphi_{b\partial(a)}u+\sum_{p=1}^n(-1)^{|b|+|b\partial||a|}\varphi_{a\partial_p(b)}(\varphi_{\zeta_{p}}\psi_{\partial}-\psi_{\zeta_{p}\partial})u
+\sum_{j=1}^m(-1)^{|a||b\partial|}\varphi_{a\rmd_j(b)}
(\varphi_{t_{j}^{-1}}\psi_{t_{j} \partial}-\psi_{ {\partial}})u\\
+&(-1)^{|a||b\partial|}\varphi_{ab}\psi_{\partial}u\in W.
\end{align*}
This shows that $W$ is invariant under $\psi_x$ for $x\in \fg={\dot A}\op{\dot \fk}.$ It is also clear that $W$ is invariant under $\varphi_a$ for $a\in {\dot A}.$ Finally, since $U$ is a $\T$-submodule of \[\Omega\sub N=M_{(\mu_0)}=\{v\in M\mid \psi_\partial(u)=\mu_0(\partial)u,~~\partial=1,\rmd_1,\ldots,\rmd_m\},\] first of all, for $\partial\in \{\rmd_i,\partial_k\mid 1\leq i\leq m, 1\leq k\leq n\}$, we have 
\[\hat\varphi_\partial U=X_{\bs{\bar r},\emptyset,\partial}U\sub U\quad \hbox{where } \bs{\bar r}=(-1,0,\ldots,0)\in\bbbz^{m+1},\] 
and second of all, for $u\in U$ and $\partial=1,\rmd_1,\ldots,\rmd_m,\partial_1,\ldots,\partial_n$, 
\[\psi_\partial\hat\varphi_1(u)=\hat\varphi_1\psi_\partial(u)\] showing that $\hat\varphi_1(u)\in \Omega.$
These altogether imply that $W$ is invariant under $\psi_x,\hat\varphi_x$ and $\varphi_a$ ($x\in\fg, a\in {\dot A}$). So, \[W=\sspan_\bbbc \{\varphi_{a}u\mid a\in {\dot A},u\in U\}=M.\]
Next recall  \(\ft'=\sspan_\bbbc\{1,\rmd_i\mid 1\leq i\leq m\}\) and for $0\leq i\leq m$  define
\begin{align*}
\gamma_i:&\ft'\longrightarrow \bbbc\\
&1\mapsto \d_{0,i}, \rmd_j\mapsto\d_{i,j}\quad (1\leq j\leq m).
\end{align*}
For $\mu\in (\ft')^*,$  $\bs{\bar r}=(r_0,\ldots,r_m)\in\bbbz^{m+1}$ and $I\sub\{1,\ldots,n\},$ we have  $\varphi_{t^{\bs{\bar r}}\zeta_I}M_{(\mu)}\sub M_{(\mu+\sum_{i=0}^mr_i\gamma_i)}.$ So,  since 
\[\sspan_\bbbc\{\varphi_{a}u\mid a\in {\dot A},u\in U\}=M,\] we have 
\begin{align*}
\sum_{\mu\in {(\ft')^*}}M_{(\mu)}=M=\sum_{I\sub\{1,\ldots,n\}}\sum_{\bs{\bar r}\in\bbbz^{m+1}}\varphi_{\zeta_I}\varphi_{t^{\bs\bar r}}U=\sum_{\bs{\bar r}\in\bbbz^{m+1}}\varphi_{t^{\bs\bar r}}(\underbrace{\sum_{I\sub\{1,\ldots,n\}}\varphi_{\zeta_I}U}_{\sub M_{(\mu_0)}=N}),
\end{align*}
we get that $N=M_{(\mu_0)}=\displaystyle{\sum_{I\sub\{1,\ldots,n\}}\varphi_{\zeta_I}U},$ as we desired.

(iv) Assume $U$ is a nonzero $\T$-submodule of $\Omega.$ Suppose that $I_1,\ldots,I_s$ are different  subsets of $\{1,\ldots,n\}$ and $I_{s_0}\neq \emptyset$ for some $1\leq s_0\leq s$.
We claim that there are no nonzero elements 
 $u_1,\ldots,u_s\in U$ such that  $u:=\varphi_{\zeta_{I_1}}u_1+\cdots+\varphi_{\zeta_{I_s}}u_s$ is an element of $\Omega$. To the contrary, assume we have such a $u.$
Let $s=1$ and $I_{1}=\{j_1<\cdots<j_p\}.$ For $1\leq r\leq p,$ set $J_r:=\{j_1<\cdots< j_{r-1}<j_{r+1}<\cdots<j_p\}$. Then, we have  
\begin{align*}
0=\psi_{\partial_{j_r}}u=&\psi_{\partial_{j_r}}\varphi_{\zeta_{I_{1}}}u_1=\psi_{\partial_{j_r}}\varphi_{\zeta_{j_1}}\cdots\varphi_{\zeta_{j_p}}u_1=(-1)^p\varphi_{\zeta_{j_1}}\cdots\varphi_{\zeta_{j_p}}\psi_{\partial_{j_r}}u_1+(-1)^{r-1}\varphi_{\zeta_{J_r}}u_1\\
=&(-1)^{r-1}\varphi_{\zeta_{J_r}}u_1.
\end{align*}
This together with an induction process on $|I_{1}|$ implies that $u_1=0$ which is a contradiction.  We then use induction on $s$ to complete the proof of the claim.

Now, if $U\neq \Omega,$ we pick  $u\in \Omega\setminus U.$ Since \[u\in \Omega\setminus U\sub N\setminus U=\sspan_\bbbc\{\varphi_{\zeta_I}u\mid I\sub\{1,\ldots,n\},u\in U\}\setminus U,\] there are nonzero elements $u_1,\ldots,u_s\in U$ and   different  subsets  $I_1,\ldots,I_s$ of $\{1,\ldots,n\}$ with  $I_{s_0}\neq \emptyset$ for some $1\leq s_0\leq s$ such that 
 $u=\varphi_{\zeta_{I_1}}u_1+\cdots+\varphi_{\zeta_{I_s}}u_s$ which is a contradiction.
 This implies that    $\Omega$ is a simple $\T$-module. 
 
For the last assertion, since $\Omega$ is a simple $\T$-module, it is a simple $\bs{S}\D$-module by  Lemma~\ref{XL1}(ii)  and so  $\bs{S}^2\D\Omega=\{0\}$ by part (ii) and  Theorem~\ref{2-simple}.    This together with Lemmas~\ref{gl} \& \ref{XL1} and (\ref{MM-30}) implies that $\Phi$ is a representation of $\frak{gl}(m+1,n)$.  By part~(ii), $\Omega$ has a weight space decomposition with respect to $\{\zeta_k\partial_k\mid 1\leq k\leq n\}$ with finite dimensional weight spaces and bounded weight multiplicities. Fix a weight space $K$ of $\Omega.$ $E_{i,i}$'s ($0\leq i\leq m$) are commuting endomorphisms on the finite dimensional vector space $K,$ so they have  a common eigenvector $v_0$.  Since $\frak{gl}(m+1,n)$ has a root space decomposition with respect to $\{E_{\a,\a}\mid 0\leq \a\leq m+n\}$ and $\Omega$ is a simple $\frak{gl}(m+1,n)$-module, the submodule  generated by $v_0$ is the entire $\Omega$ and has a weight space decomposition with respect to $\{E_{\a,\a}\mid 0\leq a\leq m+n\}$ with finite dimensional weight spaces and bounded weight multiplicities.

(v) Recall $\Omega$ and $\mu_0$ from the statement and set \[\bs{\bar \varrho}:=(\mu_0(1),\mu_0(\rmd_2),\ldots,\mu_0(\rmd_m),0,\ldots,0).\] Using Lemma~\ref{equalities}(5), it is not hard to see that the map 
\begin{align*}
\Theta:& M_{\dot A}(\Omega,\bs{\bar \varrho})\longrightarrow M\\
&a\ot\omega\mapsto\varphi_a(\omega)\quad(a\in \dot A, \omega\in \Omega)
\end{align*}
is an isomorphism.
\qed

\begin{rem}
{\rm  Suppose that $(\mathscr{B},\circ)$ is a simple  cuspidal  Harish-Chandra $A\fk$-module. Using the same  argument as in Theorem~\ref{main}, by the use of the $\dot A$-action of  ``~$\circ$~" in place of $\varphi$ and the use of $\fk$-action of ``~$\circ$~" in place of $\psi$, one can show that $(\mathscr{B},\circ)$ is a shen-Larsson module. In particular, there is a one-to-one correspondence between  simple  cuspidal  Harish-Chandra $A\fk$-modules and simple  cuspidal QP-modules over $(\dot A,\dot \fk)$. 

Next assume $V$ is a nontrivial simple cuspidal Harish-Chandra  $\fk$-module. By \cite[\S~5]{BFIK}, there is a  cuspidal Harish-Chandra $A\fk$-module $\hat V=\Bigop{\lam\in\ft^*}{}\hat V^\lam$,  where $\ft:=\sspan_\bbbc\{t_i\frac{d}{dt_i},\zeta_k\frac{\partial}{\partial\zeta_k}\mid 0\leq i\leq m, 1\leq k\leq n\}$, called the cover of $V$ and an $\fk$-module epimorphism $\pi:\hat V\longrightarrow V.$ 
We know that $\mathfrak{A}:=A\smpro U(\fk)$ has a weight space decomposition $\mathfrak{A}=\Bigop{\lam\in\ft^*}{}\mathfrak{A}^\lam$ with respect to $\ft.$
 Fix $\lam\in\ft^*$ with $\pi(\hat V^\lam)\neq \{0\}$ and note that $\hat V^\lam$ is a $\mathfrak{A}^0$-module. Suppose\footnote{This is motivated by \cite[Thm.~4.4]{LX}.} $M\sub \hat V^{\lam}$ is a minimal $\mathfrak{A}^0$-submodule of $\hat V^\lam$ with $\pi(M)\neq \{0\}$ and $M'\sub M$ is a maximal $\mathfrak{A}^0$-submodule of $M$. Then,  $\pi$ induces a  $\fk$-module epimorphism from $\mathscr{P}:=\mathfrak{A}M/\mathfrak{A}M'$ onto $V.$
The quotient  $\mathscr{P}$ is a simple cuspidal Harish-Chandra  $A\fk$-module and so, it is a Shen-Larsson $A\fk$-module as we mentioned above. So,  the nontrivial simple  cuspidal Harish-Chandra $\fk$-module $V$ is a simple  $\fk$-quotient of  a Shen-Larsson $A\fk$-module. }
\end{rem}

\end{document}